\documentclass{amsart}
\usepackage{fullpage}
\usepackage{stmaryrd}
\usepackage{amssymb}
\usepackage{mathrsfs}
\usepackage[colorlinks]{hyperref}
\usepackage{booktabs}
\usepackage{algorithm}
\usepackage{algpseudocode}

\newtheorem{theorem}{Theorem}[section]
\newtheorem{lemma}[theorem]{Lemma}

\theoremstyle{definition}
\newtheorem{definition}[theorem]{Definition}

\theoremstyle{remark}

\numberwithin{equation}{section}

\newcommand{\FF}{{\mathbb{F}}}

\newcommand{\bC}{{\mathbf{C}}}

\newcommand{\bZ}{{\mathbf{Z}}}
\newcommand{\bF}{{\mathbf{F}}}

\newcommand{\Aut}{{\operatorname{Aut}}}

\newcommand{\Sp}{{\operatorname{Sp}}}

\newcommand{\GL}{{\operatorname{GL}}}
\newcommand{\Gal}{{\operatorname{Gal}}}

\newcommand{\GF}{\mbox{GF}}

\makeatletter
\newenvironment{breakablealgorithm}
  {
   \begin{center}
     \refstepcounter{algorithm}
     \hrule height.8pt depth0pt \kern2pt
     \renewcommand{\caption}[2][\relax]{
       {\raggedright\textbf{\ALG@name~\thealgorithm} ##2\par}%
       \ifx\relax##1\relax 
         \addcontentsline{loa}{algorithm}{\protect\numberline{\thealgorithm}##2}%
       \else 
         \addcontentsline{loa}{algorithm}{\protect\numberline{\thealgorithm}##1}%
       \fi
       \kern2pt\hrule\kern2pt
     }
  }{
     \kern2pt\hrule\relax
   \end{center}
  }
\makeatother
\begin{document}

\title{On the orbits of a finite solvable primitive linear group}


\author{Yong Yang}

\author{Mengxi You}


\address{Department of Mathematics, Texas State University, 601 University Drive, San Marcos, TX 78666, USA.}
\email{yang@txstate.edu}

\address{College of Science, China Three Gorges University, Yichang, Hubei 443002, China.}
\email{mengxiyou1@hotmail.com}

\makeatletter

\makeatother

\subjclass[2000]{20C15, 20C20}
\keywords{solvable groups, primitive linear groups, orbits}
\date{}



\begin{abstract}
In this paper, we strengthen a result of Seager regarding the number of orbits of a solvable primitive linear group.
\end{abstract}

\maketitle
\maketitle
\Large
\section{Introduction} \label{sec:introduction}
Results on the orbits in finite group actions, besides being of interest in their own right, have long been known to be crucial to establish results in the character theory or the conjugacy classes of finite groups. In this paper, we strengthen a result of Seager on the action of solvable primitive linear groups. In \cite{Seager1} and \cite{Seager2}, Seager estimates the rank of a solvable primitive permutation group. One of the main results shows that, if a solvable primitive linear group is not semi-linear, then it will have at least $\frac{p^{n/2}} {12n} + 1$ orbits on the vector space, except possibly if $p^n$ is one of ten exceptional numbers. This result has many applications in counting the number of the conjugacy classes of a finite group (see for example \cite{HK2000,HK2003,HHKM2011,Keller2,KellerMoreto,Maroti}), as well as studying character degrees ~\cite{Keller1}.

In this note, we show that the list of the exceptional cases in ~\cite{Seager1} can be removed completely, and thus a non semi-linear solvable primitive linear group will always have at least $\frac {p^{n/2}} {12n} + 1$ orbits on the vector space it acts on.  This result can not only simplify the proof of several known results, but also has some potential applications. The main ingredients of the proof are some more delicate counting arguments as well as GAP ~\cite{GAP2020} and Magma calculations.


\section{Lemma and a result of Seager} \label{sec:Lemma and a result of Seager}
If $V$ is a finite vector space of dimension $n$ over $\GF(q)$, where $q$ is a prime power, we denote by $\Gamma(q^n)=\Gamma(V)$ the semilinear group of $V$, i.e.,
\[\Gamma(q^n)=\{x \mapsto ax^{\sigma}\ |\ x \in \GF(q^n), a \in \GF(q^n)^{\times}, \sigma \in \Gal(\GF(q^n)/\GF(q))\},\] and we define \[\Gamma_0(q^n)=\{x \mapsto ax\ | \ x \in \GF(q^n), a \in \GF(q^n)^{\times}\}.\]
Now we state the structure of the primitive solvable linear groups which will be used later.

\begin{definition} \label{defineEi}
Suppose that a finite solvable group $G$ acts faithfully, irreducibly and quasi-primitively on a finite vector space $V$. Let $\bF(G)$ be the Fitting subgroup of $G$ and $\bF(G)=\prod_i P_i$, $i=1, \dots, m$ where $P_i$ are normal $p_i$-subgroups of $G$ for different primes $p_i$. Let $Z_i = \Omega_1(\bZ(P_i))$. We define \[E_i=\left\{ \begin{array}{lll} \Omega_1(P_i) & \mbox{if $p_i$ is odd}; \\ \lbrack P_i,G,\cdots, G \rbrack & \mbox{if $p_i=2$ and $\lbrack P_i,G,\cdots, G \rbrack \neq 1$}; \\  Z_i & \mbox{otherwise}. \end{array} \right.\] By proper reordering we may assume that $E_i \neq Z_i$ for $i=1, \dots, s$, $0 \leq s \leq m$ and $E_i=Z_i$ for $i=s+1, \dots, m$. We define $E=\prod_{i=1}^s E_i$, $Z=\prod_{i=1}^s Z_i$ and $\overline{E}_i=E_i/Z_i$, $\bar{E}=E/Z$. Furthermore, we define $e_i=\sqrt {|E_i/Z_i|}$ for $i=1, \dots, s$ and $e=\sqrt{|E/Z|}$.
\end{definition}

\begin{theorem} \label{Strofprimitive}
Suppose that a finite solvable group $G$ acts faithfully, irreducibly and quasi-primitively on an $n$-dimensional finite vector space $V$ over finite field $\FF$ of characteristic $r$. We use the notation in Definition ~\ref{defineEi}. Then every normal abelian subgroup of $G$ is cyclic and $G$ has normal subgroups $Z \leq U \leq F \leq A \leq G$ such that,
\begin{enumerate}
\item $F=EU$ is a central product where $Z=E \cap U=\bZ(E)$ and $\bC_G(F) \leq F$;
\item $F/U \cong E/Z$ is a direct sum of completely reducible $G/F$-modules;
\item $E_i$ is an extraspecial $p_i$-group for $i=1,\dots,s$ and $e_i=p_i^{n_i}$ for some $n_i \geq 1$. Furthermore $(e_i,e_j)=1$ when $i \neq j$ and $e=e_1 \dots e_s$ divides $n$, also $\gcd(r,e)=1$;
\item $A=\bC_G(U)$ and $G/A \lesssim \Aut(U)$, $A/F$ acts faithfully on $E/Z$;
\item $A/\bC_A(E_i/Z_i) \lesssim \Sp(2n_i,p_i)$;
\item $U$ is cyclic and acts fixed point freely on $W$ where $W$ is an irreducible submodule of $V_U$;
\item $|V|=|W|^{eb}$ for some integer $b$;
\item $G/A$ is cyclic and $|G:A| \mid \dim(W)$. $G=A$ when $e=n$;
\item Let $g \in G \backslash A$, assume that $o(g)=t$ where $t$ is a prime and let $|W|=r^m$. Then $t \mid m$ and we can view the action of $g$ on $U$ as follows, $U \leq \FF_{r^{m}}^*$ and $g \in \Gal(\FF_{r^{m}}:\FF_r)$.
\end{enumerate}
\end{theorem}
\begin{proof}
This is \cite[Theorem 2.2]{YY3}.
\end{proof}

\begin{lemma} \label{OrderG}
Suppose that a finite solvable group $G$ acts faithfully,irreducibly and quasi-primitively on a finite vector space $V$. Using the notation in Theorem ~\ref{Strofprimitive}, we have $|G| \leq \dim(W) \cdot |A/F| \cdot e^2 \cdot (|W|-1)$.
\end{lemma}
\begin{proof}
By Theorem ~\ref{Strofprimitive}, $|G|= |G/A||A/F||F|$ and $|F|=|E/Z||U|$. Since $|G/A| \mid \dim(W)$, $|E/Z|=e^2$ and $|U| \mid (|W|-1)$, we have $|G| \mid \dim(W) \cdot |A/F| \cdot e^2 \cdot (|W|-1)$.
\end{proof}
\begin{lemma} \label{aforder}
Using the notation in Theorem \ref{Strofprimitive}, we have  $|A/F|$ bounds for several different values of $e$ as follows.
\begin{table}[htbp]
    \centering
    \caption{Upper bounds on $|A/F|$}
    \begin{tabular}{cc}
        \toprule  
        $e$&$|A/F| \leq$ \\
        \midrule  
        $3$&$24$ \\
        $4$&$6^2 \cdot 2$ \\
        $6$&$24 \cdot 6$ \\
        $8$&$6^4$ \\
        $9$&$24^2 \cdot 2$ \\
        $16$&$6^4 \cdot 24$\\
        \bottomrule  
    \end{tabular}
\end{table}
\end{lemma}
\begin{proof}
This is a part of \cite[Proposition 2.10]{Dey}.
\end{proof}

In ~\cite{Seager1}, Seager showed the following result.

\begin{theorem} \label{thmold}
Let $p$ be a prime and $G$ be a solvable primitive subgroup of the linear group $\GL(n,p)$. Let $r$ be the number of orbits of $G$ on the underlying vector space of $p^n$ elements. Suppose that in this action $G$ is not permutation isomorphic to a subgroup of $\Gamma(p^n)$. Then $r > \frac {p^{n/2}} {12n} + 1$ except possibly when $p^n= 17^4, 19^4, 7^6, 5^8, 7^8, 11^8, 13^8, 7^9, 3^{16}$, and $5^{16}$.
\end{theorem}

This result has several nice applications (see for example ~\cite{HK2000,HK2003,HHKM2011,Keller1,Keller2,KellerMoreto,Maroti}). While reading Seager's paper, the first author noticed that the listed exceptional cases might be superfluous. In this paper, we show that this is indeed the case. Since some of the cases are tight, we have to rely on the computational program GAP ~\cite{GAP2020} or Magma to construct the linear groups explicitly as well as calculate the exact number of orbits. The calculation is based on the methods developed in ~\cite{HoltYang}.

Our main result is the following.


\begin{theorem} \label{thmnew}
Let $p$ be a prime and let $G$ be a solvable primitive subgroup of the linear group $\GL(n,p)$. Let $r$ be the number of orbits of $G$ on the underlying vector space of $p^n$ elements. Suppose that in this action $G$ is not permutation isomorphic to a subgroup of $\Gamma(p^n)$. Then $r > \frac {p^{n/2}} {12n} + 1$.
\end{theorem}







\section{Proof of the Main result} \label{sec:Main result}
We use the notation in Theorem ~\ref{Strofprimitive}. For $\GL(n,p)$ let $n = \dim(V)=b \cdot e$. If $e$ is a prime power an $b=1$ we using these parameters and based on the methods in \cite{HoltYang} to obtain the following two algorithm.

\begin{breakablealgorithm}
\caption{}\label{extraspecial}
\begin{algorithmic}[1]
\State Construct the extraspecial group $E$ as guaranteed by Theorem \ref{Strofprimitive}
\State Construct normalizer $NE$ as subgroups of $\GL(e,p)$
\State Construct the normalizer $N$ of $NE$ in $\GL (e,p)$
\State $CongClassesN \gets [N]$
\State $lb \gets 0$
\While{$CongClassesN$ is no empty}
    \State $G$ $\gets$ Any group of $CongClassesN$
    \State $CongClassesN \gets$ Remove the group $G$ in $CongClassesN$
    \If{$G$ is  primitive}
        \If{$G$ is solvable and  irreducible}
            \If{$E \lesssim G$}
                \State $norbits$ $\gets$ the number of orbits of $G$ on $V$
            \Else
                \State Skip
            \EndIf
            \If{$lb=0$}
                \State $lb \gets norbits$
            \Else
                \If{$norbits < lb$}
                    \State $lb \gets norbits$
                \EndIf
                \State Append maximal subgroups of $G$ to $CongClassesN$
            \EndIf
        \EndIf
    \Else
        \State Skip
    \EndIf
\EndWhile\\
\Return{$lb$}
\end{algorithmic}
\end{breakablealgorithm}\par
The $lb$ obtained by Algorithm \ref{extraspecial} is the minimum number of orbits of $G$ on $V$, where $G$ is a solvable primitive subgroup of $\GL(n,p)$. \par
The following algorithm works the same as Algorithm \ref{extraspecial}, but requires less computer memory. The disadvantage of this algorithm is that it does not calculate the exact lower bound.
\begin{breakablealgorithm}
\caption{}\label{extraspecial2}
\begin{algorithmic}[1]
\State Construct the extraspecial group $E$ as guaranteed by Theorem \ref{Strofprimitive}
\State Construct normalizer $NE$ as subgroups of $\GL(e,p)$.
\State Construct the normalizer $N$ of $NE$ in $\GL (e,p)$.
\State $CongClassesN \gets [N]$
\State $GrpList \gets [\quad]$
\If{$\frac {p^{n/2}} {12n}$ is an integer}
    \State $LB=[\frac {p^{n/2}} {12n}] + 2$($[\quad]$ is integer ceiling function)
\Else
    \State $LB=[\frac {p^{n/2}} {12n}] + 1$($[\quad]$ is integer ceiling function)
\EndIf
\While{$CongClassesN$ is no empty}
    \State $G$ $\gets$ Any group of $CongClassesN$
    \State $CongClassesN \gets$ Remove the group $G$ in $CongClassesN$
    \If{$G$ is  primitive}
        \State $norbits$ $\gets$ the number of orbits of $G$ on $V$
        \If{$G$ is solvable, irreducible and $norbits < LB$}
             \If{$E \lesssim G$}
                 \State Append $G$ to $GrpList$
            \EndIf
        \Else
             \State Append maximal subgroups of $G$ to $CongClassesN$
        \EndIf
    \Else
        \State Skip
    \EndIf
\EndWhile\\
\Return{$GrpList$}
\end{algorithmic}
\end{breakablealgorithm}\par
Each group in the $GrpList$ obtained by Algorithm \ref{extraspecial2} is a primitive solvable subgroup of $\GL(n,p)$ and the number of orbits on $V$ is less than or equal to $\frac {p^{n/2}} {12n} + 1$.\par
For $\GL(6,7)$ if $e=6$ is not a prime power, we use another method to calculate.
\begin{lemma} \label{GL67}
Let $G$ be a primitive maximal solvable subgroup of $\GL(6,7)$. If $e=6$, then $G$ is conjugate in $\GL(6,7)$ to the Kronecker product $G_2 \times G_3$, where $G_2$ is a primitive solvable subgroup of $\GL(2,7)$, and $G_3$ is a primitive solvable subgroup of $\GL(3,7)$.
\end{lemma}
\begin{proof}
By \cite [Lemma 2.1(6)]{Dol08}.
\end{proof}
Based on Lemma ~\ref{GL67}, we obtain the following algorithm.
\begin{breakablealgorithm} \caption{}\label{e6}
\begin{algorithmic}[1]
\State $Group2 \gets $ All primitive solvable subgroups of $\GL(2,7)$
\State $Group2 \gets $ All primitive solvable subgroups of $\GL(3,7)$
\State $G \gets [\quad]$
\State $lb \gets 0$
\For{$G_2$ in $Group2$}
    \For{$G_3$ in $Group3$}
        \State Append the Kronecker product $G_2 \times G_3$ to $Grp$
    \EndFor
\EndFor
\For{$G_0$ in $G$}
    \If{$G_0$ is  primitive  solvable and  irreducible}
        \State $norbits$ $\gets$ the number of orbits of $G_0$ on $V$
        \If{$lb=0$}
            \State $lb \gets norbits$
        \Else
            \If{$norbits < lb$}
                \State $lb \gets norbits$
            \EndIf
        \EndIf
    \EndIf
\EndFor\\
\Return{$lb$}
\end{algorithmic}
\end{breakablealgorithm}\par
The $lb$ obtained by Algorithm \ref{GL67} is the minimum number of orbits of $G$ on $V$, where $G$ is a solvable primitive subgroup of $\GL(6,7)$.\\

We now prove our main result, Theorem ~\ref{thmnew}.
\begin{proof}
In the following we use the notation in Theorem ~\ref{Strofprimitive}. In view of the previous result, we only need to go through the exceptional cases. Since the bound was obtained using the estimate in the case $e=2$, we may assume $e>2$. In the following, we estimate $|G|$ by Lemma ~\ref{OrderG} and  $|A/F|$ by Lemma ~\ref{aforder}. For the $r$ in Theorem  ~\ref{thmnew}, we have one immediate lower bound is $\frac{|V|-1}{|G|} +1 $. This follows from the facts that at least one trivial orbit exists, that all orbits of $G$ on $V$ must partition $V$, and that the largest possible orbit size is $|G|$.\par
If this lower bound is not large enough, we can use another method to get a more precise lower bound of $r$. For each orbit $O$ of $G$, $|O|$ divides $|G|$. $G$ has at least one trivial orbit on $V$.  Let $d_{1}, d_{2}, . . . , d_{t}$ be the divisors of $|G|$ from largest to smallest. Let $n_{i}$ be a non-negative integer$(1 \leq i \leq t)$ such that $\displaystyle\sum_{i=1}^{t}n_{i}d_{i}=|V|-1$, let $N=\displaystyle\sum_{i=1}^{t}n_{i}$. Our goal is to choose the appropriate $d_{i}$ and $n_{i}$ so that $N$ is minimal. This is the change-making problem, thus we can use the method of ~\cite {Pea} to achieve our goal. When our goal is achieved, $N+1$ is the more precise lower bound of $r$, where $N$ is the theoretical minimum number of nontrivial orbits of $G$ on $V$.
\begin{enumerate}
\item  $p^n= 17^4$, $\lceil 17^2/(12 \cdot 4)\rceil+1=8$. Let $e=4$,
using Algorithm \ref{extraspecial} to calculate this case in MAGMA, we see that $G$ will have at least 19 orbits on $V$.\\

\item $p^n= 19^4$, $\lceil 19^2/(12 \cdot 4)\rceil+1=9$.

Let $e=4$, using Algorithm \ref{extraspecial} to calculate this case in MAGMA, we see that $G$ will have at least 21 orbits on $V$.
\\

\item $p^n= 7^6$, $\lceil 7^3/(12 \cdot 6)\rceil+1=6$ .

 Let $e=3$, thus $|W|=7$ and $\dim(W)=1$ or $|W|=7^2$ and $\dim(W)=2$ and $|A/F| \leq 24$. Since $|G| \leq (7^2-1) \cdot 3^2 \cdot 24 \cdot 2$. In this case, $(|V|-1)/|G|+1 > 6 $ and the result follows.\par
Let $e=6$, using Algorithm \ref{e6} to calculate this case in GAP ~\cite{GAP2020}, we see that $G$ will have at least $8$ orbits on $V$.\\

\item $p^n= 5^8$, $\lceil 5^4/(12 \cdot 8)\rceil+1=8$. Let $e=4$, thus $|W|=5$ and $\dim(W)=1$ or $|W|=5^2$ and $\dim(W)=2$ and $|A/F| \leq 72$. Since $|G| \leq (5^2-1) \cdot 2^4 \cdot 72 \cdot 2$. In this case, $(|V|-1)/|G|+1 > 8$ and the result follows.\par
Let $e=8$, using Algorithm \ref{extraspecial} to calculate this case in GAP ~\cite{GAP2020}, 
we see that $G$ will have at least $12$ orbits on $V$. \\

\item $p^n= 7^8$, $\lceil 5^4/(12 \cdot 8)\rceil+1=27$.

Let $e=4$, thus $|W|=7$ and $\dim(W)=1$ or $|W|=7^2$ and $\dim(W)=2$, $|A/F| \leq 72$ and $|G| \leq (7^2-1) \cdot 2^4 \cdot 72 \cdot 2$. In this case, $(|V|-1)/|G|+1 > 53$ and the result follows.\par
Let $e=8$, using Algorithm \ref{extraspecial} to calculate this case in GAP  ~\cite{GAP2020}, 
we see that $G$ will have at least $34$ orbits on $V$.\\

\item $p^n= 11^8$, $\lceil 11^4/(12 \cdot 8)\rceil+1=155$.

Let $e=4$, thus $|W|=11$ and $\dim(W)=1$ or $|W|=11^2$ and $\dim(W)=2$ and $|A/F| \leq 72$. Let $e=8$, thus $|W|=11$, $\dim(W)=1$ and $|A/F| \leq 6^4$. Since in all these cases, $|G| \leq 10 \cdot 2^6 \cdot 6^4$, we have $(|V|-1)/|G|+1 > 259$ and the result follows.\\

\item $p^n= 13^8$, $\lceil 13^4/(12 \cdot 8)\rceil+1=299$.
Let $e=4$, thus $|W|=13$ and $\dim(W)=1$ or $|W|=13^2$ and $\dim(W)=2$ and $|A/F| \leq 72$. Let $e=8$, thus $|W|=13$, $\dim(W)=1$ and $|A/F| \leq 6^4$. Since $|G| \leq 12 \cdot 2^6 \cdot 6^4$, we have $(|V|-1)/|G|+1 > 820$ and the result follows.\\

\item $p^n= 7^9$, $\lceil 7^4.5/(12 \cdot 9)\rceil+1=60$. Let $e=3$, thus $|W|=7$ and $\dim(W)=1$ or $|W|=7^3$ and $\dim(W)=3$ and $|A/F| \leq 24$. Let $e=9$, thus $|W|=7$, $\dim(W)=1$ and $|A/F| \leq 24^2 \cdot 2$. Since in all these cases $|G| \leq 6 \cdot 3^4 \cdot 24^2 \cdot 2$, we have $(|V|-1)/|G|+1 > 73$ and the result follows.\\

\item $p^n= 3^{16}$, $\lceil 3^8/(12 \cdot 16)\rceil+1=36$.

Let $e=4$, thus $|W|=3$ and $\dim(W)=1$, $|W|=3^2$ and $\dim(W)=2$ or $|W|=3^4$ and $\dim(W)=4$ and $|A/F| \leq 6^2 \cdot 2$. Since $|G| \leq (3^4-1) \cdot 2^4 \cdot 72 \cdot 4$. In this case, $(|V|-1)/|G|+1 > 117$ and the result follows.\par
Let $e=8$, thus $|W|=3$ and $\dim(W)=1$ or $|W|=3^2$ and $\dim(W)=2$ and $|A/F| \leq 6^4$. Since $|G| \leq (3^2-1) \cdot 2^6 \cdot 6^4 \cdot 2$ and $(|V|-1)/|G|+1  < 34$.
So we need to find a more precise lower bound using the change-making problem algorithm in Python, we see that $G$ will have at least 38 orbits on $V$.\par

Let $e=16$, using Algorithm \ref{extraspecial2} to calculate this case in the MAGMA, we see that $GrpList$ is none. Thus, $G$ will have at least 36 orbits on $V$.\\

\item $p^n= 5^{16}$,  $\lceil 5^8/(12 \cdot 16)\rceil+1=2036$.

Let $e=4$, thus $|W|=5$ and $\dim(W)=1$, $|W|=5^2$ and $\dim(W)=2$ or $|W|=5^4$ and $\dim(W)=4$ and $|A/F| \leq 6^2 \cdot 2$. Let $e=8$, thus $|W|=5$ and $\dim(W)=1$ or $|W|=5^2$ and $\dim(W)=2$ and $|A/F| \leq 6^4$. Let $e=16$, thus $|W|=5$, $\dim(W)=1$ and $|A/F| \leq 6^4 \cdot 24$. Since in all those cases, $|G| \leq 4 \cdot 2^8 \cdot 6^4 \cdot 24$, we have $(|V|-1)/|G|+1 > 4791$ and the result follows.
\end{enumerate}
\end{proof}

\noindent
\textbf{Remark}: We want to discuss how tight the $\frac {p^{n/2}} {12n} + 1$ bound is. We first notice the current bound is obtained by estimating the largest group orders under the case $e=2$ and it seems hard to find an alternative form for all the cases. Second, we have constructed a primitive solvable linear group in $\GL(16,3)$ with the following structure: $e=16$, $A/F \cong S_3 \wr S_4$, and $|W|=3$. The group action has 49 orbits while $\lceil 3^8/(12 \cdot 16)\rceil+1=36$.

\section{Acknowledgement} \label{sec:Acknowledgement}
This work was partially supported by a grant from the Simons Foundation (No. 918096, to YY). The authors would also like to thank Professor Derek Holt for his invaluable help.\\ 

\noindent \textbf{Data availability Statement:} Data sharing not applicable to this article as no datasets were generated or analysed during the current study.\\

\noindent \textbf{Competing interests:} The authors declare none.


\end{document}